\documentclass{article}

\usepackage{delarray,verbatim,enumerate,a4wide}
\usepackage{amsmath,amsthm,amstext,amsbsy,amssymb,amsfonts,amscd}
\usepackage[all]{xy}
\usepackage[french]{babel}
\usepackage[T1]{fontenc}
\usepackage[utf8]{inputenc}

\usepackage{scalerel}
\usepackage{upgreek,bbm}
\usepackage[small]{titlesec}

\usepackage{color}
\definecolor{vert}{rgb}{0.1,0.4,0.2}
\usepackage[colorlinks=true,linkcolor=blue,citecolor=vert]{hyperref}

\usepackage{calligra}
\DeclareFontShape{T1}{calligra}{m}{n}{<->s*[0.95]callig15}{}
\DeclareMathAlphabet{\mathscr}{T1}{calligra}{m}{n}

\newtheorem{Th}{Théorème}[]
\newtheorem{Lem}[Th]{Lemme}
\newtheorem{Prop}[Th]{Proposition}
\newtheorem{Cor}[Th]{Corollaire}

\newtheorem{Sco}[Th]{Scolie}
\newtheorem{Def} [Th]{Définition}

\newtheorem*{Th*}{Théorème}
\newtheorem*{Cor*}{Corollaire}
\newtheorem*{Def*}{Définition}
\newtheorem*{ThA}{Théorème A}
\newtheorem*{ThB}{Théorème B}

\def\Preuve{\noindent {\it Preuve.~}}
\def\PreuveTh{\noindent {\it Preuve du Théorème.~}}
\def\PreuveLem{\noindent {\it Preuve du Lemme.~}}
\def\Remarque{\smallskip\noindent {\it Remarque.~}}

\def\Nota{\smallskip\noindent {\it Nota.~}}

		\def\QQ{\mathbb Q}	
\def\NN{\mathbb N}	\def\ZZ{\mathbb Z}		
\def\F2{\mathbb{F}_2}	\def\Z2{\mathbb{Z}_2}		
\def\Zl{{\mathbb{Z}_\ell}} 	\def\Ql{\mathbb{Q}_\ell}		

 				\def\U{\mathcal  U}	\def\F{\mathcal  F}	\def\N{\mathcal  N}
\def\J{\mathcal  J}  		\def\C{\mathcal  C}		\def\R{\mathcal  R}	
 	\def\Pl{\mathcal  P\ell}  	\def\Cl{\mathcal  C\!\ell}	
\def\E{\mathcal  E}		\def\T{\mathcal  T}

		\def\p{{\mathfrak p}}		\def\x{{\mathfrak x}}		\def\a{{\mathfrak a}}		\def\m{{\mathfrak m}}
		\def\l{{\mathfrak l}}

\def\wi{\widetilde}				\def\Tr{\operatorname{Tr}}
	\def\div{\operatorname{div}}
	\def\deg{\operatorname{deg}}		\def\Ind{\operatorname{Ind}}
\def\Gal{\operatorname{Gal}}	\def\Log{\operatorname{Log}}		\def\Rad{\operatorname{Rad}}
\def\Ker{\operatorname{Ker}}		\def\Hom{\operatorname{Hom}}	\def\Ima{\operatorname{Im}}

\newcommand\scale[2]{\vstretch{#1}{\hstretch{#1}{#2}}}

\newcommand\si[1]{\scale{.7}{#1}}	
\newcommand\ph{{\phantom{*}}}
\newcommand\ab{{\scale{.8}{\rm ab}}}	\newcommand\lc{{\scale{.8}{\rm lc}}}	\newcommand\nr{{\scale{.8}{\rm nr}}}
\newcommand\aug{{\scale{.8}{\rm aug}}}	\newcommand\reg{{\scale{.8}{\rm reg}}}
	
\def\%{{\scale{.8}{\infty}}}

\makeatletter
\newcommand*\wt[2][0.2ex]{%
        \begingroup
        \mathchoice{\wt@helper{#1}{#2}{\displaystyle}{\textfont}}
                   {\wt@helper{#1}{#2}{\textstyle}{\textfont}}
                   {\wt@helper{#1}{#2}{\scriptstyle}{\scriptfont}}
                   {\wt@helper{#1}{#2}{\scriptscriptstyle}{\scriptscriptfont}}%
        \endgroup
        #2%
}
\newcommand*\wt@helper[4]{%
        \def\currentfont{\the#41}%
        \def\currentskewchar{\char\the\skewchar\currentfont}%
        \setbox\tw@\hbox{\currentfont$#2$\currentskewchar}%
        \dimen@ii\wd\tw@
        \setbox\tw@\hbox{\currentfont$#2${}\currentskewchar}%
        \advance\dimen@ii-\wd\tw@
        \rlap{\raisebox{-#1}{$\m@th#3\kern\dimen@ii\widetilde{\phantom{#2}}$}}%
}
\makeatother

\def\wE{\,\wt[0.1ex]{\!\mathcal E}}		\def\wU{\wt[0.2ex]{\mathcal U}}
\def\wJ{\,\wt[0.2ex]{\!\mathcal J}}	\def\wCl{\wt[0.1ex]{\mathcal C\!\ell}} \def\wDl{\wt[0.2ex]{\mathcal D\!\ell}}
	\def\wdiv{\wt[0.3ex]{\div}}			
			
\def\wF{\,\wt[0.1ex]{F}}
\begin{document}

\title{\Large\bf Annulateurs circulaires des groupes de classes logarithmiques}

\author{ Jean-François {\sc Jaulent} }
\date{}
\maketitle
\bigskip\bigskip

{\small
\noindent{\bf Résumé.} Étant donnés un corps abélien réel $F$ de groupe $G_F$ et un nombre premier impair $\ell$, nous définissons le sous-groupe circulaire $\,\wE^\circ_F$ du pro-$\ell$-groupe des unités logarithmiques $\,\wE^\ph_F$ et nous montrons que pour tout morphisme galoisien $\rho$ de $\,\wE^\ph_F$ dans $\Zl[G_F]$, l'image $\rho(\wE^\circ_F)$ annule le $\ell$-groupe des classes logarithmiques $\,\wCl_F^\ph$. Nous en déduisons une preuve de l'analogue logarithmique de la conjecture de Solomon.
}

\

{\small
\noindent{\bf Abstract.}
 Given a real abelian field $F$ with group $G_F$ and an odd prime number $\ell$, we define the circular subgroup $\,\wE^\circ_F$ of the pro-$\ell$-group of logarithmic units $\,\wE^\ph_F$ and we show that for any Galois morphism $\rho\;:\;\wE^\ph_F\to\Zl[G_F]$, the ideal $\rho(\wE^\circ_F)$ annihilates the $\ell$-group of logarithmic classes $\,\wCl_F^\ph$. We deduce from this a proof of the logarithmic version of Solomon conjecture.
}

%\tableofcontents

%\newpage
%%%%%%%%%%%%%%%%%%%%%%%%%%%%%%%%%%%%%%%%%%%%%%%%%%%%%%%%
%%%%%%%%%%%%%%%%%%%%%%%%%%%%%%%%%%%%%%%%%%%%%%%%%%%%%%%%
\section*{Introduction}
\addcontentsline{toc}{section}{Introduction}
%%%%%%%%%%%%%%%%%%%%%%%%%%%%%%%%%%%%%%%%%%%%%%%%%%%%%%%%
%%%%%%%%%%%%%%%%%%%%%%%%%%%%%%%%%%%%%%%%%%%%%%%%%%%%%%%%

%À écrire. Penser à définir $\zeta_n$ et $\approx$\smallskip
Les $\ell$-groupes de classes logarithmiques $\,\wCl_F$ introduits dans \cite{J28} sont des invariants arithmétiques canoniques attachés à un corps de nombres $F$ pour chaque nombre premier $\ell$.

Ces groupes, qui sont finis sous la conjecture de Gross-Kuz'min, donc inconditionnellement pour $F$ abélien, s'apparentent aux groupes de classes de diviseurs des corps de fonctions (ils s'obtiennent comme quotient du groupe des diviseurs de degré nul par son sous-groupe principal), sont  liés aux noyaux sauvages de la $K$-théorie et interviennent naturellement en théorie d'Iwasawa. Par exemple, la conjecture de Greenberg sur la trivialité des invariants structurels attachés aux $\ell$-groupes de classes dans la $\Zl$-extension cyclotomique $F_\%/F$ revient à postuler que le $\ell$-groupe $\,\wCl_F$ des classes logarithmiques de $F$ capitule dans $F_\%$ pour $F$ totalement réel (cf. \cite{J60,J64}).\medskip

La première ambition de ce travail est de transposer aux groupes de classes logarithmiques des corps abéliens réels les théorèmes d'annulation obtenus par Thaine, Rubin, All et al. (cf. \cite{Al1,BM,BN2,GK1,GK2,NV,Ru,So,Th}) pour les groupes de classes habituels. Procédant à la manière de Sinnott, nous introduisons pour cela un sous-groupe circulaire convenable $\,\wE_F^\circ$ du pro-$\ell$-groupe des unités logarithmiques $\,\wE_F$. C'est l'objet de la section 1.\smallskip

Ce point acquis, la section 2 suit l'approche de Rubin \cite{Ru} et de All \cite{Al1} pour établir (Th. \ref{TPC}):

\begin{ThA}
Étant donné un corps abélien réel $F$ de groupe de Galois $G_F$ et $\ell$ un premier impair, pour tout morphisme galoisien $\rho$ de $\,\wE_F$ dans $\Zl[G_F]$, l'image $\rho(\wE_F^\circ)$ du sous-groupe des unités logarithmiques circulaires annule le $\ell$-groupe des classes logarithmiques $\,\wCl_F$.
\end{ThA}

Nous en déduisons une version logarithmique d'un résultat conjecturé par Solomon (Th. \ref{Solomon}):

\begin{ThB}
Soient $F$ un corps abélien réel, $G_F$ son groupe de Galois, $\ell$ un nombre premier impair, $\l$ une place de $F$ au-dessus de $\ell$ et $\ZZ_\l$ l'anneau des entiers du complété $\l$-adique $F_\l$; puis $\vartheta$ le $\Zl[G_F]$-morphisme du pro-$\ell$-groupe $\,\wE_F$ des unités logarithmiques de $F$ dans $F_\l[G_F]$ défini par:\smallskip

\centerline{$\vartheta(\varepsilon)=\sum_{\sigma\in G_F}\Log_\l(\varepsilon^\sigma)\sigma^{\si{-1}}$.}\smallskip

\noindent Alors, pour tout élément $\a$ de l'algèbre $F_\l[G_F]$ tel qu'on ait $\a\vartheta(\wE_F)\subset\ZZ_\l[G_F]$, l'image $\a\vartheta(\wE_F^\circ)$ du pro-$\ell$-groupe des unités logarithmiques circulaires annule $\ZZ_\l\otimes_\Zl\wCl_F$.
\end{ThB}

Rappelons enfin qu'en présence des racines $\ell$-ièmes de l'unité les éléments de Stickelberger et leurs reflets permettent également de construire des annulateurs galoisiens pour les $\ell$-groupes de classes logarithmiques (cf. \cite{J65}), à la manière de Gras et Oriat (cf. \cite{Gra4,Kuc,Or1}).

%%%%%%%%%%%%%%%%%%%%%%%%%%%%%%%%%%%%%%%%%%%%%%%%%%%%%%%%%%%%
%%%%%%%%%%%%%%%%%%%%%%%%%%%%%%%%%%%%%%%%%%%%%%%%%%%%%%%%

\newpage
%%%%%%%%%%%%%%%%%%%%%%%%%%%%%%%%%%%%%%%%%%%%%%%%%%%%%%%%%%%
%%%%%%%%%%%%%%%%%%%%%%%%%%%%%%%%%%%%%%%%%%%%%%%%%%%%%%%%%%%
\section{Éléments circulaires et unités logarithmiques}
%%%%%%%%%%%%%%%%%%%%%%%%%%%%%%%%%%%%%%%%%%%%%%%%%%%%%%%%%%%
%%%%%%%%%%%%%%%%%%%%%%%%%%%%%%%%%%%%%%%%%%%%%%%%%%%%%%%%%%%
\subsection{Rappel  sur les classes et unités logarithmiques}
%\section*{Préliminaire: rappel  sur les classes et unités logarithmiques}
%\addcontentsline{toc}{section}{Préliminaire: rappel  sur les classes et unités logarithmiques}
%%%%%%%%%%%%%%%%%%%%%%%%%%%%%%%%%%%%%%%%%%%%%%%%%%%%%%%%%%%

Soient $\ell$ un nombre premier donné et $F$ un corps de nombres.  À chaque place finie  $\p$ de $F$, il est attaché dans \cite{J28} une application à valeurs dans $\Zl$, définie sur le groupe multiplicatif $F_\p^\times$ du complété de $F$ en $\p$ par la formule:%\smallskip

\centerline{$\tilde\nu_\p (x_\p)\, =\,\nu_\p (x_\p)$, pour $\p\nmid\ell$;\quad et \quad $\tilde\nu_\p (x_\p)\, =  -\frac{1}{\deg\, \p}\,\Log_\ell (N_{F_\p/\QQ_p}(x_\p))$, pour $\p|\ell$;}\smallskip

\noindent où $\Log_\ell$ désigne le logarithme d'Iwasawa et $\deg\p$ est un facteur de normalisation, dont l'expression exacte est sans importance ici, destiné à assurer que l'image de $F_\p^\times$ soit dense dans $\Zl$. Cette application induit un morphisme surjectif du {\em compactifié $\ell$-adique} du groupe multiplicatif $F_\p^\times$\smallskip

\centerline{$\R_{F_\p}=\varprojlim  F_\p^\times/F_\p^{\times\ell^n}$}\smallskip

\noindent dont le noyau, dit {\em sous-groupe des unités logarithmiques} de $\R_{F_\p}$,\smallskip

\centerline{$\wU_{F_\p}\,=\,\{x_\p\in\R_{F_\p}\,|\, \tilde\nu_\p(x_\p)=0\}$}\smallskip

\noindent  s'identifie par la Théorie $\ell$-adique locale du corps de classes (cf. \cite{J31}) au sous groupe normique de $\R_{F_\p}$ associé à la $\Zl$-extension cyclotomique $F_\p^{\,c}$ de $F_\p$.
C'est donc l'analogue du groupe\smallskip

\centerline{$\U_{F_\p}\,=\,\{x_\p\in\R_{F_\p}\,|\, \nu_\p(x_\p)=0\}$}\smallskip

\noindent des unités de de $\R_{F_\p}$, qui correspond, lui, à la $\Zl$-extension non-ramifiée de $F_\p$.
\medskip

Soit maintenant $\J _F$ le {\em $\ell$-adifié  du groupe  des idèles de $F$}, i.e. le produit $\J_F=\prod_\p ^{\si{\rm res}}\R_{F_\p}$ des compactifés  $\R _{F_\p}$ des groupes multiplicatifs des complétés $F_\p$, restreint aux familles $(x_\p)_\p$ dont presque tous les éléments tombent dans le sous-groupe unité $\,\U_F=\prod_\p \, \U _{F_\p}$.
La Théorie $\ell$-adique globale du corps de classes (cf. \cite{J31}) assure l'existence d'un isomorphisme de groupes topologiques compacts entre le $\ell$-groupe des classes d'idèles $\,\C_F$ défini comme quotient\smallskip

\centerline{$\C_F=\J_F/\R_F$}\smallskip

\noindent de $\J_F$ par son sous-groupe principal $\R_F=\Zl\otimes_\ZZ F^\times$ et le groupe de Galois $G_F=\Gal(F^{\ab}/F)$ de la pro-$\ell$-extension abélienne maximale de $F$. Dans la correspondance ainsi établie (cf. \cite{J28,J31}):\smallskip

\begin{itemize}
\item[(i)] Le groupe de normes associé à la $\Zl$-extension cyclotomique $F^c=F_{\%}$ de $F$ est le sous-groupe des idèles de degré nul: $\wJ_F=\{\x=(x_\p)_\p\in\J_F\,|\, \deg(\x)=\sum_\p\wi\nu_\p(x_\p)\deg\,\p=0\}$.\smallskip

\item[(ii)] Le groupe de normes associé à la plus grande sous-extension $F^{\lc}$ de $F^{\ab}$ qui est localement cyclotomique (i.e. complètement décomposée sur $F^c$ en chacune de ses places) est le produit $\,\wU_F\R_F$ du sous-groupe $\,\wU_F=\prod_\p \, \wU _{F_\p}$ des unités logarithmiques locales et de $\R_F$.\smallskip

\item[(iii)] En particulier, le groupe de Galois $\Gal(F^{lc}/F^c)$ s'identifie au quotient $\,\wCl_F= \wJ_F/\wU_F\R_F$, lequel peut être regardé comme quotient du groupe $\wDl _F = \wJ _F /\wU_F$ des diviseurs logarithmiques de degré nul par son sous-groupe principal $\Pl_F=\R_F\wU_F/\wU_F$, le numérateur $\,\wDl _F$  s'identifiant au sous-groupe $\wi\oplus_\p \,\Zl\,\p$ des diviseurs de degré nul de la somme formelle $\oplus_\p \,\Zl\,\p$.\smallskip

\item[(iv)] Et le noyau $\wE _F=\R_F \cap\, \wU_F$ du morphisme $\wdiv:\;x\mapsto \sum_\p\wi\nu_\p(x)\,\p$ de $\R_F$ dans $\wDl_F$ est le sous-groupe des normes cyclotomiques (locales comme globales) de $\R_F$.
\end{itemize}
\smallskip

\begin{Def*}
Nous disons que $\,\wCl_F= \wJ_F/\wU_F\R_F \simeq \wDl_F/\Pl_F$ est le $\ell$-groupe des classes logarithmiques (de degré nul) du corps $F$ et que $\,\wE_F$ est le pro-$\ell$-groupe des unités logarithmiques globales.\smallskip

Le quotient $\,\wCl^*_F=\J_F/\wU_F\R_F$, qui s'identifie non canoniquement à la somme directe de $\wCl_F$ et de $\Zl$, est, par convention, le pro-$\ell$-groupe  des classes logarithmiques de degré arbitraire.
\end{Def*}

Comme expliqué dans \cite{J28}, la {\em conjecture de Gross-Kuz'min} (pour le corps $F$ et le premier $\ell$) revient à postuler la finitude du (pro)-$\ell$-groupe $\,\wCl_F$ ou, de façon équivalente, que le $\Zl$-rang du pro-$\ell$-groupe des unités logarithmiques $\,\wE_k$ est le somme $r_{\si{F}}+c_{\si{F}}$ des nombres de places réelles et complexes de $F$. Elle est toujours vérifiée dès lors que $F$ est abélien.\smallskip

Enfin, du point de vue de la théorie d'Iwasawa, le groupe $\,\wCl_F$ s'interprète comme le quotient des genres ${}^\Gamma\T_F$, relativement au groupe procyclique $\Gamma=\Gal(F_\%/F)$, du {\em module de Kuz'min-Tate}\smallskip

\centerline{$\T_F = \varprojlim \,\Cl'_{F_n}$}\smallskip

\noindent limite projective des $\ell$-groupes de $\ell$-classes d'idéaux attachés aux étages finis $K_n$ de la tour $K_\%/K$.

\newpage

%%%%%%%%%%%%%%%%%%%%%%%%%%%%%%%%%%%%%%%%%%%%%%%%%%%%%%%%%%%
\subsection{Éléments cyclotomiques logarithmiques}
%%%%%%%%%%%%%%%%%%%%%%%%%%%%%%%%%%%%%%%%%%%%%%%%%%%%%%%%%%%

Notons $\ell$ un nombre premier arbitraire et $(\zeta_n)_{n>1}$ un système cohérent de racines primitives de l'unité, en ce sens qu'on ait $\zeta_m^{m/n}=\zeta_n$ pour $n|m$, par exemple en posant $\zeta_n=\exp 2i\pi/n$.

\begin{Prop}
Soit $F$ un corps abélien de conducteur $f>1$ et de groupe $G_F=\Gal(F/\QQ)$.
\begin{itemize}\smallskip
\item Pour $\ell \mid f$, l'élément $\eta_F^\ph=N_{\QQ[\zeta_f]/F}(1-\zeta_f^\ph)$ est une unité logarithmique: $\eta_F^\ph\in\wE_F^\ph$.\smallskip
\item Pour $\ell \nmid f$, l'intersection du $\Zl[G_F]$ module multiplicatif engendré par $\eta_F$ avec le pro-$\ell$-groupe $\wE_F$ des unités logarithmiques de $F$ contient le $\Zl[G_F]$-module engendré par l'élément

\centerline{$\wi\eta_F^\ph = \eta_F^{1-(\frac{F}{\ell})_\ph^{\si{-1}}}$.}\smallskip
\end{itemize}

En d'autres termes, on a: $\;\eta_F^{\,\Zl[G_F]} \subset \;\wE_F$ dans le premier cas; $\wi\eta_F^{\,\Zl[G_F]} \subset \;\wE_F$ dans le second.
\end{Prop}

\Preuve Il est bien connu (cf. e.g. \cite{Gra4}, \S 4.2) que les $\eta_F^\ph$ satisfont les identités normiques:
\begin{equation}\label{IdNorm}
N_{F/K}(\eta_F^\ph)=\eta_K^{\prod_p\,\big( 1-\big(\frac{K}{p}\big)^{\si{-1}}\big)}
\end{equation}

\noindent où, pour toute sous-extension $K$ de $F$, le produit fait intervenir les symboles d'Artin attachés aux premiers $p$ qui se ramifient dans $F/\QQ$ mais non dans $K/\QQ$. 

Rappelons en outre que les $\eta_F^\ph$ sont des $p$-unités qui ne sont pas unités lorsque $f$ est  $p$-primaire pour un premier $p$; des unités sinon.
Cela étant, distinguons les cas:\smallskip

\begin{itemize}

\item Pour $\ell\mid f$, le cas $f=p^r$ avec $p\ne\ell$ étant exclu, $\eta_F^\ph$ est toujours une $\ell$-unité et il vient donc:\smallskip

\centerline{$\wi v_\p(\eta_F^\ph)=v_\p(\eta_F^\ph)=0$, pour chaque $\p\nmid\ell$,}\smallskip

puisqu'aux places étrangères à $\ell$ les valuations logarithmique $\wi v_\p$ et ordinaire $v_\p$ coïncident.

Aux places au-dessus de $\ell$, il vient, en revanche:\smallskip

\centerline{$\wi v_\l(\eta_F^\ph)=\frac{1}{\deg\l}\Log_\ell(N_{F_\l/\Ql}(\eta_F^\ph))$.}\smallskip

Introduisons donc le sous-corps de décomposition $F_\circ$ de $\ell$ dans $F/\QQ$. Par hypothèse le premier $\ell$ se ramifie alors dans $F$ mais non dans $F_\circ$ et l'on a en outre $\big(\frac{F_\circ}{\ell}\big)=1$, de sorte que la formule normique plus haut nous donne immédiatement:\smallskip

\centerline{$\wi v_\l(\eta_F^\ph)=\frac{1}{\deg\l}\Log_\ell(N_{F/F_\circ}(\eta_F^\ph))=\frac{1}{\deg\l}\Log_\ell 1=0$,}\smallskip

sauf dans le cas $\ell$-primaire $f=\ell^r$, où l'on a directement:\smallskip

\centerline{$\wi v_\l(\eta_F^\ph)=\frac{1}{\deg\l}\Log_\ell(N_{F/\QQ}(\eta_F^\ph))=\frac{1}{\deg\l}\Log_\ell \ell=0$,}\smallskip

En fin de compte, il vient bien $\wi v_\l(\eta_F^\ph)=0$ pour $\l\mid\ell$; et $\eta_F^\ph$ est une unité logarithmique.\smallskip

\item Pour $\ell\nmid f$, l'élément $\eta_F^\ph$ est toujours une unité (et donc une unité logarithmique aux places en dehors de $\ell$), sauf dans le cas primaire $f=p^r$ où c'est une $p$-unité qui n'est pas une unité.
Dans cette dernière situation, si $\p$ est alors l'unique place de $F$ au-dessus de $p$, l'égalité $\wi v_\p(\eta_F^\alpha)= v_\p(\eta_F^\alpha)$, pour tout $\alpha$ dans $\Zl[G_F]$ donne l'équivalence:\smallskip

\centerline{$\wi v_\p(\eta_F^\alpha)=0\;\Leftrightarrow\;\alpha\in \ZZ_\ell^{\si{\rm aug}}[G_F]$, idéal d'augmentation de l'algèbre $\Zl[G_F]$.}
\smallskip

Reste dans tous les cas à évaluer les valuations logarithmiques aux places $\l$ au-dessus de $\ell$.\smallskip

La formule $\wi v_\l(\eta_F^\alpha)=\frac{1}{\deg\l}\Log_\ell(N_{F_\l/\Ql}(\eta_F^\alpha))$
donne l'équivalence:\smallskip

\centerline{$\wi v_\l(\eta_F^\alpha)=0\;\Leftrightarrow\;N_{F/F_\circ}(\eta_F^\alpha)\in\ell^\Zl\;\Leftrightarrow\;N_{F/F_\circ}(\eta_F^\alpha)=1$,}\smallskip

où $F_\circ$ désigne, comme précédemment, le sous-corps de décomposition de $\ell$. Et cette dernière condition est évidemment remplie lorsque $\alpha$ est contenu dans l'idéal d'augmentation du sous-groupe de décomposition $D_\ell$ de $\ell$; autrement dit lorsque c'est un multiple de $ 1-(\frac{F}{\ell})^{\si{-1}}_\ph$.
\end{itemize}

\begin{Sco}
Sous les mêmes hypothèses, les éléments $\eta_F^\ph$ pour $\ell\mid f$ (et $\wi\eta_F^\ph$ pour $\ell\nmid f$) sont en fait des normes universelles dans la $\Zl$-tour cyclotomique $F_\%$ de $F$, i.e. des éléments de l'intersection $\N_F=\bigcap_{n\in\NN}N_{F_n/F}(\wE_{F_n})$ des groupes de normes logarithmiques attachés aux étages finis de  $F_\%$.
\end{Sco}

\Preuve C'est une conséquence immédiate des identités normiques citées plus haut. On a, en effet: $N_{F_n/F}(\eta_{F_n}^\ph)=\eta_F^\ph$, dans le premier cas; $N_{F_n/F}(\eta_{F_n}^\ph)=\wi\eta_F^\ph$, dans le second.

\newpage
%%%%%%%%%%%%%%%%%%%%%%%%%%%%%%%%%%%%%%%%%%%%%%%%%%%%%%%%%%%
\subsection{Unités logarithmiques circulaires}
%%%%%%%%%%%%%%%%%%%%%%%%%%%%%%%%%%%%%%%%%%%%%%%%%%%%%%%%%%%

Supposons maintenant $\ell$ impair et prenons encore $F$ abélien réel de conducteur $f=f_K$ et de groupe de Galois $G_F=\Gal(F/\QQ)$. Notons enfin $F_\%=\bigcup_{n\in\NN}F_n$ sa $\Zl$-extension cyclotomique. \smallskip

Le $\ell$-adifié $\R^\circ_F$ du groupe des éléments circulaires ({\em à la Sinnott}) de $F$ est  le $\Zl[G_F]$-module engendré dans $\R_F=\Zl\otimes F^\times$ par les images des éléments $\eta_K^\ph=N_{\QQ[\zeta_{f_{\si{K}}}]/K}(1-\zeta_{f_{\si{K}}})$, pour $K\subset F$.\smallskip

Il pourrait donc paraître naturel de définir le pro-$\ell$-groupe des unités logarithmiques circulaires du corps abélien réel $F$ comme l'intersection $\,\wE_F^\ph\cap\,\R_F^\circ$ de $\R_F^\circ$ avec le pro-$\ell$-groupe $\,\wE_F$ des unités logarithmiques. Malheureusement les éléments ainsi obtenus ne satisfont pas clairement la propriété spéciale de Rubin qui joue un rôle essentiel dans la preuve du Théorème \ref{TPC} (cf. Lem. \ref{LC} infra). C'est pourquoi il est préférable de procéder comme suit:

\begin{Def}
Soient $\ell$ un nombre premier impair, $F$ un corps abélien réel, $f=f_F$ son conducteur, $G_F=\Gal(F/\QQ)$ son groupe de Galois et $F_\%=\bigcup_{n\in\NN}F_n$ sa $\Zl$-extension cyclotomique.\par

Pour $K \subset F$ de conducteur $f_K$, notons $\eta_K^\ph=N_{\QQ[\zeta_{f_{\si{K}}}]/K}(1-\zeta_{f_{\si{K}}})$ et $\wi\eta_K^\ph = \eta_K^{1-(\frac{K}{\ell})_\ph^{\si{-1}}}$.\smallskip

Le pro-$\ell$-groupe des unités logarithmiques circulaires (à la Sinnott) de $F$ est le $\Zl[G_F]$-module $\,\wE^\circ_F$ engendré conjointement par les images des $\eta^\ph_K$ pour $\ell\mid f_K$ et des $\wi\eta^\ph_K$ pour $\ell\nmid f_K$. 
\end{Def}

\Nota Les identités $N_{F_n/F}(\eta_{K_n}^\ph)=\eta_K^\ph$ pour $\ell\mid f_K$ et  $N_{F_n/F}(\eta_{K_n}^\ph)=\wi\eta_K^\ph$ pour $\ell\nmid f_K$ montrent que les unités logarithmiques circulaires sont en particulier des normes universelles; i.e. que l'on a:\smallskip

\centerline{$\,\wE^\circ_K\subset\N_F=\bigcap_{n\in\NN}N_{F_n/F}(\wE_{F_n})$.}

\begin{Th}\label{Car}
Soient $F$ un corps abélien réel et $F^\circ$ le sous-corps de décomposition de $\ell$. Notons $G_F$ le groupe de Galois $\Gal(F/\QQ)$ et $D_F$ le sous-groupe $\Gal(F/F^\circ)$.\smallskip

\begin{itemize}
\item[(i)] Lorsque $F$ possède plus d'une place au-dessus de $\ell$, i.e. pour $F^\circ\ne\QQ$, le caractère $\wi\chi_F^\circ$ du $\Zl[G_F]$-module des unités logarithmiques circulaires $\,\wE^\circ_F$ est alors l'induit à $G_F$ du caractère d'augmentation du sous-groupe de décomposition $D_F$:\smallskip

\centerline{$\wi\chi_F^\circ= \Ind_{D_F}^{G_F}\,\chi^\aug_{D_F}$.}\smallskip

\item[(ii)] Le même résultat $\wi\chi_F^\circ=\chi^\aug_{G_F}$ vaut encore avec $D_F=G_F$ lorsque le corps $F$ admet une unique place au-dessus de $\ell$, sauf si $F$ contient un sous-corps $K$ de conducteur $\ell$-primaire, auquel cas $\wi\chi_F^\circ=\chi^\reg_{G_F}$ est le caractère régulier.
\end{itemize}
\end{Th}

La preuve de ce résultat repose sur le Lemme suivant:

\begin{Lem}\label{déc}
Pour tout corps abélien réel $F\ne\QQ$ dans lequel le nombre premier impair $\ell$ se décompose complètement, le pro-$\ell$-groupe des unités logarithmiques circulaires  est trivial: $\,\wE^\circ_F=1$.
\end{Lem}

\PreuveLem En l'absence de ramification en $\ell$ dans $F/\QQ$, les $\eta^\ph_K$ pour $K\subset F$ sont des unités aux places au-dessus de $\ell$. On a donc: $\,\wE_F^\circ\subset\,\E^\ph_F\cap\,\wE_F^\ph$. Puis, sous l'hypothèse de complète décomposition: $\,\E^\ph_F\cap\,\wE_F^\ph=\mu_F^\ph=1$ (cf. \cite{Grt}, p. 218, l. 11--13, ou \cite{J64}, Prop.6).\smallskip

\PreuveTh Regardons d'abord $(i)$. Le Lemme donne $N_{F\!/\!F^\circ}(\wE^\circ_F)\subset\wE^\circ_{F^\circ}=\mu^\ph_{F^\circ}=1$. De $\,\wE^\circ_F \subset\,\wE^\ph_F$, qui est ici un $\Zl[G_F]$-module de caractère $\chi_{G_F}^\reg$ régulier,on tire donc: $\wi\chi_F^\circ \le \Ind_{D_F}^{G_F}\,\chi^\aug_{D_F}$.\par

Or, par construction  le groupe $\,\wE^\circ_F$ contient l'image $\,\E^\circ_F{}^{I_{D_F}}$du pro-$\ell$-groupe des unités circulaires $\,\E^\circ_F$ par l'idéal d'augmentation $I_{D_F}=\sum_{\sigma\in D_F}\Zl[G_F](\sigma-1)$.
Comme $\,\E^\circ_F$, qui est d'indice fini dans $\,\E_F^\ph$, est un $\Zl[G_F]$-module de caractère $\chi_{G_F}^\aug$, on conclut: $\wi\chi_F^\circ \ge \Ind_{D_F}^{G_F}\,\chi^\aug_{D_F}$. D'où le résultat.\smallskip

Examinons maintenant $(ii)$, qui correspond à $F^\circ=\QQ$. Les formules normiques pour les unités logarithmiques circulaires donnent $N_{K/\QQ}(\eta_{\si{K}}^\ph)=1$ dans le $\ell$-adifié $\R_\QQ$, pour $K\subset F$ sauf si $f_{\si{K}}$ est $\ell$-primaire, auquel cas il vient $N_{K/\QQ}(\eta_{\si{K}}^\ph)=\ell$.\smallskip

Lorsque ce cas est exclu, on a donc $\wi\chi_F^\circ\le\chi^\aug_{G_F}$ et on conclut comme précédemment à l'égalité.\smallskip

En revanche, si $F$ contient un sous-corps $K$ de conducteur $\ell$-primaire, $\,\wE^\circ_F$ contient l'image de $\ell$ et $\wi\chi_F^\circ$ contient le caractère unité. Comme on a $\wi\chi_F^\circ\ge\chi^\aug_{G_F}$, il vient alors $\wi\chi_F^\circ=\chi^\reg_{G_F}$.

\newpage
%%%%%%%%%%%%%%%%%%%%%%%%%%%%%%%%%%%%%%%%%%%%%%%%%%%%%%%%%%%
%%%%%%%%%%%%%%%%%%%%%%%%%%%%%%%%%%%%%%%%%%%%%%%%%%%%%%%%%%%
\section{Annulateurs circulaires des groupes de classes logarithmiques}
%%%%%%%%%%%%%%%%%%%%%%%%%%%%%%%%%%%%%%%%%%%%%%%%%%%%%%%%%%%
%%%%%%%%%%%%%%%%%%%%%%%%%%%%%%%%%%%%%%%%%%%%%%%%%%%%%%%%%%%
\subsection{Application du Théorème de \v Cebotarev}
%%%%%%%%%%%%%%%%%%%%%%%%%%%%%%%%%%%%%%%%%%%%%%%%%%%%%%%%%%%

%Reprenons dans le cadre logarithmique considéré ici les arguments de Rubin (cf. \cite{Ru}, \S5). \smallskip

Prenons toujours $\ell$ impair et considérons un corps abélien réel $F$ de conducteur $f=f_F>1$. Écrivons $G_F=\Gal(F/\QQ)$ son groupe de Galois. Notons $\,\wCl_F$ le $\ell$-groupe des classes logarithmiques (de degré nul) attaché à $F$ et $\,\wCl_F^*$ le groupe des classes logarithmiques sans condition de degré.

Prenons $m$ assez grand pour avoir ${\ell^m}\wCl_F=0$, de sorte que $\,\wCl_F$ puisse être regardé canoniquement comme un sous-groupe du quotient d'exposant $\ell^m$ du pro-$\ell$-groupe des classes logarithmiques de degré arbitraire: $\,\wCl_F^\ph \subset {}^{\ell^{\si{m}}}\!\wCl^*_F=\wCl^*_F/\ell^m\wCl_F^*$.

Désignons par $\wF$ la plus grande extension abélienne de $F$ qui est d'exposant $\ell^m$ et localement cyclotomique, de sorte que nous avons par la Théorie du corps de classes: $\Gal(\wF/F)\simeq {}^{\ell^{\si{m}}}\!\wCl^*_F$. Et notons $ F_{m_{\si{0}}}=\wF \cap F_\%$ le sous-corps de $\wF$ fixé par $\,\wCl_F$.\smallskip

Donnons-nous enfin un morphisme galoisien $\bar\rho$ du pro-$\ell$-groupe $\,\wE_F$ des unités logarithmiques de $F$ vers l'algèbre $\ZZ/\ell^m\ZZ [G_F]$; et considérons les extensions emboîtées:\smallskip

\centerline{ $F_\zeta=F[\zeta_{\ell^{\si{m}}}] \subset F_\rho=F_\zeta[\sqrt[\ell^{\si{m}}]{\Ker {\bar\rho}_\ph^\ph}] \subset F_\varepsilon = F_\zeta[\sqrt[\ell^{\si{m}}]{\wE_F^\ph}]$,}\smallskip

Observons pour cela que tout élément $\mathfrak x$ du tensorisé $\R_F=\Zl\otimes_\ZZ F^\times$ peut être représenté par un élément $x$ de $F^\times$ modulo une puissance $\ell^m$-ième de $\R_F$, disons $\mathfrak x=x\mathfrak y^{\ell^{\si{m}}}$, et que cet élément $x$ est  unique modulo $F^{\times\,\ell^{\si{m}}}$; ce qui permet de définir sans ambiguïté $F_\zeta [\sqrt[\ell^{\si{m}}]{\mathfrak x}]$ comme étant $F_\zeta [\sqrt[\ell^{\si{m}}]{x}]$.

Observons en outre que les trois extensions $F_\zeta$, $F_\rho$ et $F_\varepsilon$ sont linéairement disjointes de $\wF$ sur $ F_{m_{\si{0}}}$, de sorte que leurs composés $\wF_\zeta$, $\wF_\rho$ et $\wF_\varepsilon$ avec $\wF$ donnent lieu au schéma galoisien:

\begin{displaymath}
\xymatrix{ 
&& \wF \ar@{-}_{\wCl_F}[d] \ar@{-}[r] & \wF_\zeta \ar@{-}[d] \ar@{-}[r] & \wF_\rho \ar@{-}[d] \ar@{-}[r] & \wF_\varepsilon \ar@{-}[d] &\\
&& F_{m_{\si{0}}}  \ar@{-}[d] \ar@{-}[r]\ar@{-}[r] & F_\zeta \ar@{-}^{G_{F_{\si{\zeta}}}=\Gal(F_\zeta/\QQ)}[ddl] \ar@{-}[r] &  F_\rho \ar@{-}[r]  & F_\varepsilon\\ 
&& F \ar@{-}_{G_F}[d] &  & & & &\\
&& \QQ & & & & &
}
\end{displaymath}

Par construction, le radical kummérien de $F_\varepsilon/F_\rho$ est: $\Rad(F_\varepsilon/F_\rho) = \wE_F/\Ker\bar\rho \simeq \Ima\bar\rho$. C'est un sous-module de $\ZZ/\ell^m\ZZ[G_F]$ ; et le groupe de Galois $\Gal(F_\varepsilon/F_\rho) \simeq \Hom_{G_{F_{\si{\zeta}}}}(\Rad(F_\varepsilon/F_\rho), \mu_{\ell^{\si{m}}})$, qui s'identifie donc à un quotient de  $\ZZ/\ell^m\ZZ[G_F]$, est, de ce fait,  $\Zl[G_{F_{\si{\zeta}}}]$-monogène.\smallskip

Soient alors $\sigma^\ph_{F_\varepsilon}\! \in \Gal(F_\varepsilon/F_\rho)$ un $\Zl[G_{F_{\si{\zeta}}}]$-générateur de $\Gal(\wF_\varepsilon/\wF_\rho)$ et $\sigma_{\wF}\! \in \Gal(\wF/F_{m_{\si{0}}})$ provenant d'une classe donnée arbitraire $[\mathfrak c]$ de $\,\wCl_F$ regardée dans $ {}^{\ell^{\si{m}}}\!\wCl^*_F$.\smallskip

Définissons enfin $\sigma$ dans $\Gal(\wF_\varepsilon /F_\rho)$ en imposant $\sigma_{|_{F_\varepsilon}}=\sigma^\ph_{F_\varepsilon}$ et $\sigma_{|_{\wF}}=\sigma^\ph_{\wF}$. Cela étant:

\begin{Prop}\label{Cebotarev}
Il existe une infinité de nombres premiers impairs $p\ne\ell$, complètement décomposés dans $F_\rho/\QQ$, tels que l'image de l'une des places de $\wF_\varepsilon$ au-dessus de $p$ par l'opérateur de Frobenius coïncide avec $\sigma$. Sont alors satisfaites en particulier les propriétés suivantes:
\begin{itemize}
\item[(i)] $p$ est complètement décomposé dans $F[\zeta_{\ell^{\si{m}}}]$, donc vérifie la congruence: $p\equiv 1 [\,{\rm mod}\;\ell^m]$.
\item[(ii)] La classe $[\p]$ de l'un des premiers $\p$ de $F$ au-dessus de $p$ dans ${}^{\ell^{\si{m}}}\!\wCl^*_F$ coïncide avec $[\mathfrak c]$.
\item[(iii)] Les Frobenius dans $F_\varepsilon$ des places au-dessus de $p$ engendrent $\Gal(F_\varepsilon/F_\rho)$.
\end{itemize}
\end{Prop}

\Preuve L'existence est une conséquence immédiate du théorème de densité de \v Cebotarev. De plus,
\begin{itemize}
\item[(i)] la congruence résulte de la condition de complète décomposition dans $F_\rho$, donc dans $\QQ[\zeta_{\ell^{\si{m}}}]$;
\item[(ii)] provient du fait que $\sigma^\ph_{\wF}$ est induit par l'image canonique $[\mathfrak c]\,\Cl^{*\,\ell^{\si{m}}}_F$ de $[\mathfrak c]$ dans $ {}^{\ell^{\si{m}}}\!\wCl^*_F$;
\item[(iii)] enfin, résulte du fait que les conjugués de  $\sigma_{|_{F_\varepsilon}}$ par $G_{F_{\si{\zeta}}}$ engendrent $\Gal(F_\varepsilon/F_\rho)$.
\end{itemize}

\newpage
%%%%%%%%%%%%%%%%%%%%%%%%%%%%%%%%%%%%%%%%%%%%%%%%%%%%%%%%%%%
\subsection{Lemmes d'annulation logarithmiques}
%%%%%%%%%%%%%%%%%%%%%%%%%%%%%%%%%%%%%%%%%%%%%%%%%%%%%%%%%%%

Conservons les mêmes notations, supposons fixé le premier $\p$ de $F$ au-dessus de $p\equiv 1 [\,{\rm mod}\;\ell^m]$ donné par la Proposition \ref{Cebotarev}  et considérons  la restriction à $\,\wE_F$ du morphisme de semi-localisation $s_p$ induit par le plongement de $F$ dans le produit de ses complétés $F_{\p^{\si{\sigma}}}$ aux places au-dessus de $p$. 

Par construction, pour chacune des places $\p^\sigma$ de $F$ au-dessus de $p$, le complété $F_{\p^{\si{\sigma}}}$ s'identifie à $\QQ_p$ et le $\ell$-sous-groupe de Sylow $\mu_{\p^{\si{\sigma}}}^\ph$ du groupe $F_{\p^{\si{\sigma}}}^\times$ est d'ordre $\ell^{m_{\si{p}}}$ avec $m_p=v_\ell(p-1)\ge m$.\par

En particulier le quotient ${}^{\ell^{\si{m}}}\!\!\mu^\ph_{F_p}=\mu^\ph_{F_q}/\mu_{F_p}^{\,\ell^{\si{m}}}$ est un $\ZZ/\ell^m\ZZ[G_F]$-module libre de dimension 1.

\begin{Lem}\label{factorisation}
Sous les hypothèses de la Proposition, le morphisme galoisien $\bar\rho\;:\;\wE_F\to \ZZ/\ell^m\ZZ[G_F]$ s'écrit $\bar\rho=\bar\rho_p\circ \bar s_p$, où $\bar s_p\;:\;\wE_F \to {}^{\ell^{\si{m}}}\!\!\mu^\ph_{F_p}$ est induite par l'application de semi-localisation et $\bar\rho_p$ est un morphisme galoisien de ${}^{\ell^{\si{m}}}\!\!\mu^\ph_{F_p}$ vers $\ZZ/\ell^m\ZZ[G_F]$.
\end{Lem}

\Preuve Le noyau de $s_p$ dans $\,\wE_F$ est formé des unités logarithmiques qui sont localement puissances $\ell^m$-ièmes aux places au-dessus de $p$, i.e. des $\varepsilon\in\wE_F$ pour lesquelles les extensions $F[\sqrt[\ell^{\si{m}}]{\varepsilon}]/F$ sont complètement décomposées aux places $\p^\sigma$ (ou, ce qui revient au même, pour lesquelles l'extension $F_\zeta[\sqrt[\ell^{\si{m}}]{\varepsilon}]/F_\zeta$ est complètement décomposée aux places au-dessus de $p$). Du fait du choix de $\sigma^\ph_{F_\varepsilon}$, cela revient à exiger que les conjugués de $\alpha^\ph_{F_\varepsilon}$ agissent trivialement sur $F_\zeta[\sqrt[\ell^{\si{m}}]{\varepsilon}]$; autrement dit, que $\varepsilon$ soit une puissance $\ell^m$-ième dans $F_\rho$, i.e. le produit d'un élément de $\Ker\bar\rho$ et d'une puissance $\ell^m$-ième  dans $F_\zeta$ et finalement dans $F$, puisque, $\ell$ étant impair, les  éléments de $F$ qui sont des puissances $\ell^m$-ièmes dans $F_\zeta$ sont {\em déjà} des puissances $\ell^m$-ièmes dans $F$ (cf. e.g. \cite{Ru}, Lem. 5.7). En fin de compte, il suit: $\Ker\bar s_p=\Ker\bar\rho$, ce qui assure, par $\ZZ/\ell^m\ZZ[G_F]$-injectivité de ${}^{\ell^{\si{m}}}\!\!\mu^\ph_{F_p}$, l'existence du morphisme factorisant $\bar\rho_q$.\medskip

Soit maintenant $\zeta_p$ une racine primitive $p$-ième de l'unité. L'extension $F[\zeta_p]/F$  possède une unique sous-extension $F'$ de degré $\ell^m$, laquelle est cyclique, totalement ramifiée au-dessus de $p$ et non-ramifiée en dehors. En particulier, si  $\p'=\p_{\si{F'}}$ désigne l'unique place de $F'$ au-dessus de $\p=\p_{\si{F}}$ et $F'_\p$ son complété en la place $\p'$, on a d'une part $\p_{\si{F}}^\ph=\p_{\si{F'}}^{\ell^{\si{m}}}$; et, d'autre part, l'égalité entre $\ell$-groupes de racines locales de l'unité: $\mu^\ph_{F^{\si'}_{\p}}=\mu^\ph_{F^\ph_\p}$.\smallskip
 
 Le lemme qui suit peut être regardé comme l'analogue logarithmique du Th.5.1 de \cite{Ru}.

\begin{Lem}\label{annulation}
Avec les notations ci-dessus, soient $\,\wE_{F'\!/F}=\{\eta^\ph_{F'}\in\wE_{F'}\;|\; N_{F'\!/F}(\eta^\ph_{F'})=1\}$ le groupe des unités logarithmiques relatives de l'extension $F'\!/F$ et $s_p$ l'homomorphisme de semi-localisation à valeurs dans la somme directe $\mu^\ph_{F'_p}=\mu^\ph_{F^\ph_p}=\oplus_{\p|p}\;\mu^\ph_{F_{\p}}\simeq \mu^\ph_{\QQ_p}\otimes\Zl[G_F] \simeq (\ZZ/\ell^{m_{\si{p}}}\ZZ)[G_F]$.\smallskip

Tout $\alpha\in\Zl[G_F]$ qui annule le quotient $\mu^\ph_{F_p}/s_p(\wE_{F'\!/F})\mu^{\ell^{\si{m}}}_{F_p}$ annule la classe de $\p$ dans $^{\ell^{\si{m}}}\!\wCl_F^*$.
\end{Lem}

\Preuve Écrivons $F'_{\p'}=F^\ph_\p[\sqrt[\ell^{\si{m}}]{\pi^\ph_\p}\,]$ pour une uniformisante $\pi_\p\in F_\p$; et $x_{\si{F'}}\in\R_{F'}$ un relèvement de $(\sqrt[\ell^{\si{m}}]{\pi^\ph_\p},1,\cdots,1)\in\prod_{\p_{\si{F'}}^\ph| p}\,\R_{F'_{p^\ph_{\si{F'}}}}$; notons enfin $\delta$ un générateur de $\Delta=\Gal(F'\!/F)\simeq\Gal(F'_{\p'}/F^\ph_\p)$.\par
 Par construction, nous avons $s_p(x_{\si{F'}}^\ph)^{\delta-1}=s_p(x_{\si{F'}}^{\delta-1}))=(\zeta_\p,1,\cdots,1)$ pour une racine $\ell^m$-ième de l'unité $\zeta_\p\in F_\p$; puis, par hypothèse, $s_p(x_{\si{F}}^\alpha)^{\delta-1}=s_p(x_{\si{F}}^{\delta-1}))^\alpha=s_p(\eta^\ph_{\si{F'}})$, pour un $\eta^\ph_{\si{F'}}\!\in\wE_{F'\!/F}$. Le Théorème 90 de Hilbert nous permet alors d'écrire $\eta^\ph_{\si{F'}}=y_{\si{F'}}^{\delta-1}$ pour un $y^\ph_{\si{F'}}\!\in\R_{F'}$, qui engendre donc un diviseur logarithmique ambige. Il suit (en notations additives):\smallskip
 
 \centerline{$\wdiv y^\ph_{\si{F'}}= \sum_{\sigma\in G_F}\,\alpha'_\sigma\p_{\si{F'}}^\sigma\;+\;\a'_{\si{F}}$;\qquad puis: \qquad $\wdiv N_{F'\!/F}(y^\ph_{\si{F'}})= \sum_{\sigma\in G_F}\,\alpha'_\sigma\p_{\si{F}}^\sigma\;+\;\ell^m\a'_{\si{F}}$}\smallskip

\noindent de sorte que l'élément $\alpha'=\sum\alpha'_\sigma\sigma\in\Zl[G_F]$ annule bien la classe de $\p=\p_{\si{F}}$ dans $\,\wCl_F^*/\ell^m\wCl_F^*$.\smallskip

Reste à vérifier que $\alpha'$ et $\alpha$ sont dans la même classe modulo $\ell^m\Zl[G_F]$. Pour cela, observons que l'identité $s_p(x_{\si{F'}}^\alpha/y_{\si{F'}}^\ph )^{(\delta-1)}=1$ nous donne $s_p(x_{\si{F'}}^\alpha)=s_p(y_{\si{F'}}^\ph z_{\si{F}}^\ph)$, pour un $z_{\si{F}}^\ph$ convenable de $\R_F$. Écrivant alors $\alpha=\sum\alpha_\sigma\sigma$ et prenant les diviseurs logarithmiques respectifs des deux membres de l'identité, nous restreignant enfin aux seules composantes au-dessus de $p$, nous obtenons ainsi:\smallskip

\centerline{$\sum_{\sigma\in G_F}\,\alpha_\sigma\p_{\si{F'}}^\sigma \; =\; \sum_{\sigma\in G_F}\,\alpha'_\sigma\p_{\si{F'}}^\sigma \;+ \a_F $}
\smallskip

\noindent pour un diviseur logarithmique $\a_F$ de $F$, donc finalement, comme attendu:\smallskip

\centerline{$\sum_{\sigma\in G_F}\,\alpha_\sigma\sigma \; \equiv\; \sum_{\sigma\in G_F}\,\alpha'_\sigma\sigma \; [\,{\rm mod}\; \ell^m\Zl[G_F]\,]$,}\smallskip

\noindent puisque les diviseurs logarithmiques au-dessus de $p$ sont totalement ramifiés dans $F'\!/F$.

\newpage
%%%%%%%%%%%%%%%%%%%%%%%%%%%%%%%%%%%%%%%%%%%%%%%%%%%%%%%%%%%
\subsection{Annulation des classes logarithmiques réelles}
%%%%%%%%%%%%%%%%%%%%%%%%%%%%%%%%%%%%%%%%%%%%%%%%%%%%%%%%%%%

Nous pouvons maintenant énoncer le Théorème principal sur les annulateurs circulaires.

\begin{Th}\label{TPC}
Étant donné un corps abélien réel $F$ de groupe de Galois $G_F$ et $\ell$ un premier impair, pour tout morphisme galoisien $\rho$ de $\,\wE_F$ dans $\Zl[G_F]$, l'image $\rho(\wE_F^\circ)$ du sous-groupe des unités logarithmiques circulaires spéciales annule le $\ell$-groupe des classes logarithmiques $\,\wCl_F$.
\end{Th}

\Preuve Choisissons $m$ assez grand pour que $\ell^m$ annule le $\ell$-groupe $\,\wCl_F$ des classes logarithmiques; prenons un élément $\alpha$ dans $\rho(\wE_F^\circ)$; et partons d'une classe $[\mathfrak c ]$ de $\,\wCl_F$ regardée dans $^{\ell^{\si{m}}}\!\wCl_F^*$.\par

Notons $\bar\rho$ la réduction de $\rho$ modulo $\ell^m$, à valeurs dans $\ZZ/\ell^m\ZZ[G_F]$, et faisons choix d'un premier impair $p\ne\ell$ satisfaisant les conditions de la Proposition \ref{Cebotarev}. Par construction $p\equiv 1 [\,{\rm mod}\;\ell^m]$ est complètement décomposé dans $F$ et l'on a $[\mathfrak c]=[\p]$ dans $^{\ell^{\si{m}}}\!\wCl_F^*$ pour l'un des premiers $\p$ de $F$ au-dessus de $p$.\par

Tout revient donc à montrer que l'on a: $\bar\alpha[\p]=0$ dans $^{\ell^{\si{m}}}\!\wCl_F^*$, où $\bar\alpha\in\Ima\bar\rho$ est la réduction de $\alpha$ modulo $\ell^m$. Or, par le Lemme \ref{factorisation}, $\bar\rho$ se factorise via l'application de semi-localisation $\bar s_p$. Et:

\begin{Lem}\label{LC}
Pour toute unité logarithmique circulaire $\eta \in \wE_F^\circ$, il existe une unité logarithmique circulaire relative $\eta' \in\wE^\ph_{F'\!/F}$ telle qu'on ait: $s_p(\eta)=s_p(\eta')$.
\end{Lem}

Ainsi $\bar\alpha$, qui est l'image par $\bar\rho$ d'une unité logarithmique circulaire $\eta$, provient d'une unité logarithmique circulaire relative $\eta' \in\wE^\ph_{F'\!/F}$; et il résulte alors du Lemme \ref{annulation} qu'il annule la classe de $\p$ dans $^{\ell^{\si{m}}}\!\wCl_F^*$. D'où le résultat annoncé.\medskip

\noindent{\it Preuve du Lemme.}  Il s'agit de vérifier que les unités logarithmiques circulaires sont {\em mutatis mutandis} ce que Rubin appelle des unités spéciales dans \cite{Ru}. Reprenons pour cela dans le cadre logarithmique les calculs de All (cf. \cite{Al1}, \S3): prenons $K\subset F$, notons $f=f_{\si{K}}$ son conducteur et partons d'une unité logarithmique circulaire $\varepsilon=\eta_{\si{K}}$ (pour $\ell\mid f$) ou $\varepsilon=\wi\eta_{\si{K}}$ (pour  $\ell\nmid f)$; notons $\QQ'$ l'unique sous-corps de degré $\ell^m$ de $\QQ[\zeta_p]$ et $K'=K\QQ'$; identifions $G_F=\Gal(F/\QQ)$ à $G'_F=\Gal(F[\zeta_p]/\QQ[\zeta_p])$; et raisonnons modulo $\m'=\prod_{\p'|p}\,\p'$.\par

De $\zeta_p\equiv 1$, nous tirons: $\eta^\ph_{\si{\QQ[\zeta_f]}} = 1-\zeta_f \equiv 1-\zeta_f\zeta_p  = (1-\zeta_{pf})^\rho = \eta^\rho _{\si{\QQ[\zeta_{pf}]}}$ pour un $\rho$ de $G'_F$; et finlement: $\eta^\ph_{\si{K}} = N_{\QQ[\zeta_f]/K}(1-\zeta_f) \equiv N_{\QQ[\zeta_{pf}]/K'}(1-\zeta_{pf})^\rho= \eta_{\si{K'}}^\rho$; et, par conséquent: $\wi\eta^\ph_{\si{K}}\equiv  \wi\eta_{\si{K'}}^\rho$.\par

Posant alors $\varepsilon'=\eta_{\si{K'}}^{\rho}$  (pour $\ell\mid f$) ou $\varepsilon'=\wi\eta_{\si{K'}}^{\rho}$  (pour $\ell\nmid f$), nous obtenons bien $\varepsilon \equiv \varepsilon'$, i.e. $s_p(\varepsilon) = s_p(\varepsilon')$, puisque $\varepsilon$ est une unité aux places au-dessus de $p$, et $N_{F'\!/F}(\varepsilon')=N_{K'\!/K}(\varepsilon')=1$, puisque $p$ est ramifié dans $K'$ mais complètement décomposé dans $K$.
 \medskip

 Le résultat obtenu ci-dessus pour les classes logarithmiques peut naturellement être mis en parallèle avec celui de Rubin \cite{Ru} sur les classes d'idéaux tel que présenté par All (cf. \cite{Al1}, \S3). Désignons pour cela par $\,\Cl_F^{\,\si{0}}$ le sous-groupe du $\ell$-groupe des classes d'idéaux $\,\Cl_F^\ph$ engendré par les idéaux de degré nul, de sorte que l'on a $\,\Cl_F^{\,\si{0}}\simeq\Gal(F_\%F^\nr/F_\%)$, où $F^\nr$ désigne le $\ell$-corps de classes de Hilbert de $F$. Il vient alors:
 
 \begin{Cor}\label{CP}
 Soient $\ell$ un premier impair, $F$ abélien réel de groupe $G_F$ et $\rho$ un morphisme galoisien du $\ell$-adifié $\,\E'_F=\Zl\otimes_\ZZ E'_F$ du groupe des $\ell$-unités de $F$ dans l'algèbre $\Zl[G_F]$. Alors:\smallskip
 \begin{itemize}
 \item[(i)] L'image $\rho(\E_F^\circ)$ du pro-$\ell$-groupe construit sur les unités circulaires annule $\,\Cl_F^{\,\si{0}}$.\smallskip

  \item[(ii)] L'image $\rho(\wE_F^\circ)$ du pro-$\ell$-groupe des unités logarithmiques circulaires annule $\,\wCl_F^\ph$.
  \end{itemize}
 \end{Cor}

L'assertion $(i)$ n'est autre qu'une réécriture $\ell$-adique du résultat initial de Rubin; l'assertion $(ii)$ provient directement du Théorème \ref{TPC}. Les classes logarithmiques sont conventionnellement de degré nul (sauf mention explicite du contraire), mais non les classes au sens ordinaire; de ce fait l'introduction du sous-groupe $\,\Cl_F^{\,\si{0}}$ renforce le parallélisme des résultats.

\Remarque Lorsque le corps $F$ possède plusieurs places au-dessus de $\ell$, i.e. lorsque le sous-corps de décomposition $F^\circ$ de $\ell$ n'est pas le corps des rationnels $\QQ$, le Théorème \ref{Car} affirme en particulier que $\,\wE_F^\circ$ est contenu dans le noyau de la norme $N_{F/F^\circ}$. Le Théorème \ref{TPC} ne donne donc de ce fait aucune information directe sur le groupe $\,\wCl^\ph_{F^\circ}$.\smallskip

\newpage
%%%%%%%%%%%%%%%%%%%%%%%%%%%%%%%%%%%%%%%%%%%%%%%%%%%%%%%%%%%
%%%%%%%%%%%%%%%%%%%%%%%%%%%%%%%%%%%%%%%%%%%%%%%%%%%%%%%%%%%
\section*{Appendice: lien avec la conjecture de Solomon}
\addcontentsline{toc}{section}{Appendice: lien avec la conjecture de Solomon}
%%%%%%%%%%%%%%%%%%%%%%%%%%%%%%%%%%%%%%%%%%%%%%%%%%%%%%%%%%%
%%%%%%%%%%%%%%%%%%%%%%%%%%%%%%%%%%%%%%%%%%%%%%%%%%%%%%%%%%%
%%%%%%%%%%%%%%%%%%%%%%%%%%%%%%%%%%%%%%%%

Supposons toujours $F$ abélien réel et $\ell$ impair, mais fixons maintenant l'une $\l$ des places au-dessus de $\ell$. Solomon a conjecturé dans \cite{So} que si $\ell$ ne se ramifie pas dans $F$ l'élément\smallskip

\centerline{$\vartheta^{\si{\,\rm Sol}}_F=\frac{1}{\ell}\,\underset{\sigma\in G_{\si{F}}}{\sum}\,\Log_\l(\eta_F^\sigma)\sigma^{\si{-1}}$}\smallskip

\noindent  annule le tensorisé $\ZZ_\l\otimes_\Zl\Cl_F$ du $\ell$-groupe des classes d'idéaux de $F$.\par

Conséquence du Théorème principal de Mazur-Wiles \cite{MW} dans le cas semi-simple $\ell\nmid[F:\QQ]$, ce résultat a été prouvé par Belliard et Nguyen Quang Do dans \cite{BN1} pour $\ell$ décomposé, et sans restriction par All (cf. \cite{Al1}, Th. 1.1) sous une forme plus générale qu'on peut réécrire comme suit:

\begin{Th*}[All]
Soient $F$ un corps abélien réel, $G_F$ son groupe de Galois, $\ell$ un nombre premier impair, $\l$ une place de $F$ au-dessus de $\ell$ et $\ZZ_\l$ l'anneau des entiers du complété $\l$-adique $F_\l$; puis $\vartheta$ le $\Zl[G_F]$-morphisme du $\ell$-adifié $\,\E_F=\Zl\otimes_\ZZ E_F$ du groupe des unités dans $F_\l[G_F]$ défini par:\smallskip

\centerline{$\vartheta(\varepsilon)=\sum_{\sigma\in G_F}\Log_\l(\varepsilon^\sigma)\sigma^{\si{-1}}$.}\smallskip

\noindent Alors, pour tout  $\a$ de $F_\l[G_F]$ tel qu'on ait $\a\vartheta(\E_F)\subset\ZZ_\l[G_F]$, l'image $\a\vartheta(\E_F^\circ)$du pro-$\ell$-groupe des unités circulaires annule le tensorisé $\ZZ_\l\otimes_\Zl\Cl_F^{\,\si{0}}$ du $\ell$-groupe des classes d'idéaux de degré nul.
\end{Th*}

\Nota Dans l'isomorphisme du corps de classes $\,\Cl_F\simeq\Gal(F^\nr/F)$, où $F^\nr$ désigne la $\ell$-extension abélienne non ramifiée de $F$ (i.e. son $\ell$-corps de classes de Hilbert), le sous-groupe $\,\Cl_F^{\,\si{0}}$ des classes de degré nul correspond à $\Gal(F^\nr/(F^\nr\cap F_\%))$, où $F_\%$ est la $\Zl$-extension cyclotomique de $F$.\medskip

En parfaite analogie avec ce résultat, nous pouvons énoncer en termes logarithmiques:

\begin{Th}\label{Solomon}
Soient $F$ un corps abélien réel, $G_F$ son groupe de Galois, $\ell$ un nombre premier impair, $\l$ une place de $F$ au-dessus de $\ell$ et $\ZZ_\l$ l'anneau des entiers du complété $\l$-adique $F_\l$; puis $\vartheta$ le $\Zl[G_F]$-morphisme du pro-$\ell$-groupe $\,\wE_F$ des unités logarithmiques de $F$ dans $F_\l[G_F]$ défini par:\smallskip

\centerline{$\vartheta(\varepsilon)=\sum_{\sigma\in G_F}\Log_\l(\varepsilon^\sigma)\sigma^{\si{-1}}$.}\smallskip

Alors, pour tout élément $\a$ de l'algèbre $F_\l[G_F]$ tel qu'on ait $\a\vartheta(\wE_F)\subset\ZZ_\l[G_F]$, l'image $\a\vartheta(\wE_F^\circ)$du pro-$\ell$-groupe des unités logarithmiques circulaires annule $\ZZ_\l\otimes_\Zl\wCl_F$.

\end{Th}

\Preuve Elle est strictement identique à celle donnée dans \cite{Al1}, \S3. Rappelons-en  brièvement l'argumentation: donnons-nous une $\Zl$-base $(\mathfrak v_1,\cdots,\mathfrak v_d)$ de $\ZZ_\l$ et notons $(\mathfrak v_1^*,\cdots,\mathfrak v_d^*)$ la base duale de la codifférente. Partons d'un élément $\a=\sum a_\sigma\sigma^{\si{-1}}$ de $F_\l[G_F]$ et posons $L_\a(\varepsilon)=\sum_{\sigma\in G_{\si{F}}}a_{\sigma^{\si{-1}}}\Log_\l(\varepsilon^\sigma)$. Nous obtenons $\a\vartheta(\varepsilon)=\sum L_\a(\varepsilon^\sigma)\sigma^{\si{-1}}$, avec $L_\a(\varepsilon^\sigma)\in\ZZ_\l$ par hypothèse.\smallskip

Écrivons maintenant $L_\a(\varepsilon^\sigma)=\sum_{i=1}^d \Tr(\mathfrak v_i^*L_\a(\varepsilon^\sigma))\mathfrak v_i$ la décomposition de $L_\a(\varepsilon^\sigma)$ dans $\ZZ_\l$. Nous obtenons: $\a\vartheta(\varepsilon)=\sum_{\sigma\in G_{\si{F}}}\big(\sum_{i=1}^d  \Tr(\mathfrak v_i^*L_\a(\varepsilon^\sigma))\mathfrak v_i\big)\sigma^{\si{-1}} = \sum_{i=1}^d \big( \sum_{\sigma\in G_{\si{F}}} \Tr(\mathfrak v_i^*L_\a(\varepsilon^\sigma))\sigma^{\si{-1}}\big)\mathfrak v_i$, i.e.\smallskip

\centerline{$\a\vartheta(\varepsilon)=\sum_{i=1}^d \vartheta_{\mathfrak v_i^*}(\varepsilon) v_i$ \quad avec \quad $\vartheta_{\mathfrak v_i^*}(\varepsilon) = \sum_{\sigma\in G_{\si{F}}} \Tr(\mathfrak v_i^*L_\a(\varepsilon^\sigma))\sigma^{\si{-1}}$;}\smallskip

\noindent et l'application $\vartheta_{\mathfrak v_i^*}\;:\;\varepsilon\mapsto \sum_{\sigma\in G_{\si{F}}} \Tr(\mathfrak v_i^*L_\a(\varepsilon^\sigma))\sigma^{\si{-1}}$ est un $\Zl[G_F]$-morphisme de $\,\wE_F$ dans $\Zl[G_F]$.\smallskip

En fin de compte, si $\mathfrak x$ est un élément de $\ZZ_\l$ et $[\mathfrak c ]$ une classe de $\,\wCl_F$, il vient:\smallskip

\centerline{$\a\vartheta(\varepsilon).(\mathfrak x\otimes [ \mathfrak c ]) =\sum_{i=1}^d  \mathfrak x\mathfrak  v_i \otimes \vartheta_{\mathfrak v_i^*}(\varepsilon)[\mathfrak c]$.}\smallskip

\noindent Or, on a $\vartheta_{\mathfrak v_i^*}(\varepsilon)[\mathfrak c]=0$, si $\varepsilon$ est une unité logarithmique circulaire, en vertu du Théorème \ref{TPC}; d'où le résultat annoncé.\medskip

\Remarque Prenant $\varepsilon=\eta_F$ et $\a=\frac{1}{\ell}$, on obtient l'élément de Solomon $\frac{1}{\ell}\,\sum_{\sigma\in G_F}\Log_\l(\eta_F^\sigma)\sigma^{\si{-1}}$, lequel annule donc $\ZZ_\l\otimes_\Zl\Cl_F^{\,\si{0}}$, comme établi dans \cite{Al1}, puisque $\eta_F$ est bien une unité circulaire.\par
Si la place $\ell$ se ramifie dans $F$, $\eta_F$ est aussi une unité logarithmique circulaire. L'élément de Solomon annule alors le groupe logarithmique $\ZZ_\l\otimes_\Zl\wCl_F$. C'est en particulier le cas, dès que $F$ contient le sous-corps réel $\QQ[\zeta_\ell+\bar\zeta_\ell]$ du corps cyclotomique $\QQ[\zeta_\ell]$.

\newpage
%%%%%%%%%%%%%%%%%%%%%%%%%%%%%%%%%%%%%%%%
%Références
%%%%%%%%%%%%%%%%%%%%%%%%%%%%%%%%%%%%%%%%
\def\refname{\normalsize{\sc  Références}}

\medskip\noindent
{\small
\begin{tabular}{l}
Institut de Mathématiques de Bordeaux \\
Université de {\sc Bordeaux} \& CNRS \\
351 cours de la libération\\
F-33405 {\sc Talence} Cedex\\
courriel : Jean-Francois.Jaulent@math.u-bordeaux.fr\\
{\footnotesize \url{https://www.math.u-bordeaux.fr/~jjaulent/}}
\end{tabular}
}

\end{document}